\documentclass[11pt, reqno]{amsart}
\usepackage[cp1251]{inputenc}
\usepackage[english]{babel}
\usepackage{amssymb}
\usepackage{amsfonts}
\usepackage{amsmath}
\usepackage{mathtools}
\usepackage{amsthm}
\usepackage{xcolor}
\usepackage{hyperref}
\usepackage{amsaddr}
\usepackage{enumitem}
\usepackage{graphicx}
\usepackage{comment}

\usepackage[noabbrev]{cleveref}
\usepackage[left=1in,right=1in,top=1in,bottom=1in]{geometry}

\usepackage{array,tabularx,tabulary,booktabs}
\usepackage{longtable}
\usepackage{multirow}
\usepackage{wrapfig}
\usepackage{subcaption}
\usepackage[table]{xcolor}
\usepackage{soul}

\newtheorem{theorem}{Theorem}

\newtheorem{problem}{Problem}

\newenvironment{nouppercase}{%
  \renewcommand{\uppercasenonmath}[1]{}}{}

\begin{document}

\title{\large On non-canonical maximum cliques in the block graphs of orthogonal arrays}

\author[T.~Feng, S.~Goryainov, E.~V.~Konstantinova, and L.~Zhang]{Tao Feng$^1$, Sergey Goryainov$^2$, 
Elena V.\ Konstantinova$^{3,4}$, Letian Zhang$^5$}


\address{$^1$Beijing Jiaotong University, Beijing, P.R. China}
\address{$^2$School of Mathematical Sciences, Hebei International Joint Research Center for Mathematics and Interdisciplinary Science, Hebei Key Laboratory of Computational Mathematics and Applications, Hebei Workstation for Foreign Academicians, Hebei Normal University, Shijiazhuang 050024, P.R. China}
\address{$^3$Sobolev Institute of Mathematics, Ak. Koptyug av. 4, Novosibirsk, 630090, Russia}
\address{$^4$Novosibirsk State University, Pirogova str. 2, Novosibirsk, 630090, Russia}
\address{$^5$School of Mathematical Sciences, Hebei Normal University, Shijiazhuang 050024, P.R. China}

\email{tfeng@bjtu.edu.cn}
\email{sergey.goryainov3@gmail.com}
\email{e\_konsta@math.nsc.ru}
\email{zhangletian9@gmail.com}



\date{October 22, 2026}

\begin{abstract}
In this note, we prove an Erd\H{o}s-Ko-Rado characterisation of maximum intersecting families of columns of certain orthogonal arrays. This solves a problem of Chris Godsil and Karen Meagher and implies that, in all orthogonal arrays OA$(m,(m-1)^2)$, every non-canonical maximum intersecting family of columns necessarily has the structure of an affine plane. 
\end{abstract}


\begingroup
\begin{nouppercase}
\let\MakeUppercase\relax
\maketitle
\end{nouppercase}
\endgroup

\section{Introduction}
In their book \cite[Chapter 5]{GM15}, Godsil and Meagher proposed the studying of maximum intersecting families in classical combinatorial objects: orthogonal arrays (mutually orthogonal Latin squares) and Steiner 2-designs or, equivalently, studying maximum cliques in the associated graphs. 
An analogue of the Erd\H{o}s-Ko-Rado (EKR) theorem usually has two parts: an upper bound for the size of an intersecting family and a characterisation of the intersecting families whose size meets this bound. An important tool for establishing EKR type results is the Hoffman's ratio bound on the size of an independent set in a regular graph (see, for example, \cite[Chapter 2]{GM15} and \cite{H21}). This bound involves the degree and smallest eigenvalue of the graph. As the graphs associated with orthogonal arrays and 2-designs are strongly regular (see \cite[Chapter 10]{GR01}, \cite[Chapter 5]{GM15}, \cite{BV22}), their eigenvalues are easy to determine, and thus the Hoffman's bound applied to these graphs implies the bound in the EKR theorem for orthogonal arrays and 2-designs (the tightness can be easily shown by the existence of canonical intersecting families). Thus, to establish an EKR theorem for orthogonal arrays and Steiner 2-designs it suffices to characterise the maximum cliques in the associated graphs.
The case of Steiner of Steiner 2-designs was recently addressed in \cite{AGKK26}. 

In this note we completely resolve the following problem by Chris Godsil and Karen Meagher.

\begin{problem}[{\cite[Problem 16.4.2]{GM15}}]\label{prob:Main}
Assume that OA$(m,(m-1)^2)$ is an orthogonal array and its orthogonal array graph has non-canonical cliques of size $(m-1)^2$. Do these non-canonical cliques form subarrays that are isomorphic to orthogonal arrays with entries from $\{1,\ldots,m-1\}$?
\end{problem}

Our main result is the following theorem.

\begin{theorem}\label{thm:MainTheorem}
Assume that OA$(m,(m-1)^2)$ is an orthogonal array and its orthogonal array graph has non-canonical cliques of size $(m-1)^2$. Then all non-canonical cliques in the orthogonal array graph form complete subarrays of type OA$(m,(m-1)^2)$, that is, have the structure of affine planes.    
\end{theorem}

The paper is organised as follows. In Section \ref{sec:prelim}
we give preliminary definitions and results. In Section \ref{sec:proof} we give a proof of Theorem \ref{thm:MainTheorem}.

\section{Preliminaries}\label{sec:prelim}
We consider undirected graphs without loops and multiple edges.
A $k$-regular graph $\Gamma$ with $v$ vertices is called \emph{strongly regular} with parameters $(v,k,\lambda, \mu)$
if any two adjacent vertices in $\Gamma$ have precisely $\lambda$ common neighbours and any two distinct non-adjacent vertices in $\Gamma$ have precisely $\mu$ common neighbours.
A strongly regular graph is called \emph{imprimitive} if it is a disjoint union of cliques of the same size or a complete 
multipartite graph with parts of the same size. A strongly regular graph is called \emph{primitive} if it is not imprimitive. 
It is well known that a primitive strongly regular graph has exactly three distinct eigenvalues: the principal eigenvalue $k$ and
non-principal eigenvalues $\theta$ and $\tau$, where $\tau < 0 < \theta < k$. 
On the other hand, a connected regular graph with exactly three distinct 
eigenvalues is strongly regular (see \cite[Lemma 10.2.1]{GR01}).

A \emph{clique} in a graph $X$ is a subset of vertices of $X$ such that any two of them are adjacent.
Let $X$ be a primitive strongly regular graph $X$ with degree $k$ 
and smallest eigenvalue $\tau$. The following upper bound for the size of a clique in $X$, known as Delsarte bound 
(see \cite[Section 3.3.2]{D73}), can be derived from the Hoffman bound and implies that the size of a clique in $X$ is at most $1 - \frac{k}{\tau}$.
A clique in a strongly regular graph whose size meets the Delsarte bound is called a \emph{Delsarte clique}.

An {\em orthogonal array} $OA(n,m)$ is an $m \times n^2$ array with entries from an $n$-element set $T$
with the property that the columns of any $2 \times n^2$ subarray consist of all $n^2$ possible pairs. 
It is well known \cite[Section 10.4]{GR01} that the existence of an orthogonal array $OA(n,m)$ is equivalent to the existence of $m-2$ mutually orthogonal Latin squares of order $n$. It is also well known that for every orthogonal array $OA(n,m)$, the inequality $m \le n+1$ holds. An orthogonal array $OA(n,n+1)$ is called \emph{complete}. Further, it is well known that the existence of a complete orthogonal array $OA(n,n+1)$ is equivalent to the existence of a finite affine plane of order $n$. Finally, the existence of an finite affine plane of order $n$ is equivalent \cite{BJL99} to the existence of a finite projective plane of order $n$. Note that the question of for which positive integers there exists a projective plane of order $n$ is a major open problem in finite geometry. If $n$ is a prime power, then the existence is known. 

Two columns of an orthogonal array are called \emph{intersecting} if there exists a row in which they have the same entry (it follows from the definition that for any two columns, there is at most one such a row).
The {\em block graph of an orthogonal array} $OA(n,m)$, denoted $X_{OA(n,m)}$,  is defined to be the graph
whose vertices are columns of the orthogonal array, where two columns are adjacent if they are intersecting. Let $S_{r,i}$ be the set of columns of $OA(n,m)$ that have the entry $i$ in row $r$. These sets are cliques, and since each element of the $n$-element set $T$ occurs
exactly $n$ times in each row, the size of $S_{r,i}$ is $n$ for all $i$ and $r$. These cliques
are called the \emph{canonical cliques} in the block graph $X_{OA(m,n)}$. A simple combinatorial argument shows that the block graph of an orthogonal array is strongly regular (see \cite[Theorem 5.5.1]{GM15}). The Delsarte bound applied to the block graph of an orthogonal array $OA(n,m)$ gives the value $n$, and the existence of the canonical cliques (of size $n$) shows the tightness of the Delsarte bound for the block graphs of orthogonal arrays. 

If their book Godsil and Meagher gave a proof \cite[Corollary 5.5.3]{GM15} of the fact that if $m < \sqrt{n}+1$, then the block graph of an orthogonal array $OA(n,m)$ has only canonical maximum cliques. Thus, if
\begin{equation}\label{ineq:CharacterisationFails}
m\ge \sqrt{n}+1,    
\end{equation}
then the characterisation may fail. Godsil and Meagher posed the following problem.

\begin{problem}[{\cite[Problem 16.4.1]{GM15}}] 
Find a characterisation of the orthogonal arrays, based only
on the parameters of the array, for which all of the maximum cliques in the
orthogonal array graph are the canonical cliques.
\end{problem}

Let $m-1$ be a prime power. Then there exist a complete orthogonal array $OA(m-1,m)$. Godsil and Meagher point out that MacNeish construction \cite[Section 6.4.2]{S04} can produce an orthogonal array $OA(n,m)$ where $n = (m-1)^2$, where these parameters $n$ and $m$ correspond to equality case in the inequality (\ref{ineq:CharacterisationFails}). It is interesting that the original array $OA(m-1,m)$ corresponds to a non-canonical maximum clique in the block graph of the resulting array $OA((m-1)^2,m)$. Godsil and Meagher thus posed the following Problem \ref{prob:Main}.

Finally, the most general problem due to Godsil and Meagher \cite[Section 16.4]{GM15} is to determine all the maximum cliques in the block graph for any orthogonal array.

\section{Proof of the main theorem}\label{sec:proof}

Let $\mathcal F$ be a maximum-sized intersecting family of columns in an OA$(m,(m-1)^2)$. Then by \cite[Theorem 5.5.2]{GM15}, $|\mathcal F|=(m-1)^2$. For convenience, we call each column in $\mathcal F$ a \emph{block} (this term originates from the equivalence between orthogonal arrays and transversal designs; consequently, each element (called a \emph{point}) in a block is regarded as a pair $(i,a)$, where $i$ corresponds to the $i$-th row of the orthogonal array, and $a$ is an entry of the array).

By the proof of \cite[Theorem 5.5.2]{GM15}, for any block $F\in \mathcal F$ and any point $x$ in $F$, the number of blocks in $\mathcal F$ containing $x$ is $m-1$. 

Let $V=\bigcup_{F\in \mathcal F} F$ be the set of all points appearing in at least one member of $\mathcal F$, and let $v=|V|$. Then by double-counting the total number of occurrences of points, we have
$$v(m-1)=|\mathcal F|m,$$
and hence $v=m(m-1)$. It remains to show that $\mathcal F$ forms a sub-orthogonal array OA$(m,(m-1))$ of the given OA$(m,(m-1)^2)$.

Consider the dual set system $D=(P,\mathcal L)$, where $P=\mathcal F$ and $\mathcal{L}=\{L_p: p\in V\}$, with $L_p=\{B\in \mathcal F : p \in B\}$. Then $D$ forms an $((m-1)^2, m-1, 1)$-BIBD. Since every $(n^2,n,1)$-BIBD with $n\geq 2$ is an affine plane of order $n$ (see \cite[Chapter I, Exercise 2.21(b)]{BJL99}), the lines in $\mathcal{L}$ are partitioned into $m$ parallel classes, each consisting of $m-1$ lines. Note that each line $L_p$ corresponds to an original point $ p $; thus this partition divides the original point set $V$ into $m$ subsets $V_1,V_2,\dots,V_m$, each of size $m-1$.

Take any original block $F \in \mathcal F$; it corresponds to a point in $D$. In $D$, every point lies on exactly one line from each parallel class. Hence, $F$ lies on $m$ lines $L_{p_1}, \dots, L_{p_m}$, where these $p_i$ come from distinct $V_i$. This means $p_i \in F$, so $F$ contains exactly one point from each $V_i$.

Label the points in each $V_i$ as $1,2, \dots, m-1$. Then each block $F$ corresponds to an $m$-tuple $(a_1,a_2,\dots,a_m)$, where $a_i$ is the label of the point from $V_i$ contained in $F$.
Since any two blocks intersect in exactly one point, for any $i \neq j$, each pair of labels $(a_i, a_j)$ uniquely determines a block. Therefore, all $(m-1)^2$ $m$-tuples form an orthogonal array OA$(m,m-1)$.

\section*{Declaration of A.I.\ use}

All the content of this note is due to the authors. No A.I. tools were used to obtain the proof.

\end{document}